\documentclass{agtart_a}
\pdfoutput=1

\title{Quasiflats with holes in reductive groups}
\author{Kevin Wortman}
\givenname{Kevin}
\surname{Wortman}
\address{Department of Mathematics\\
Yale University\\
10 Hillhouse Ave\\\newline
PO Box 208283\\
New Haven CT 06520-8283\\USA} 
\email{kevin.wortman@yale.edu}



\received{18 November 2004}
\revised{}
\accepted{11 August 2005}
\proposed{}
\seconded{}
\published{24 February 2006}
\publishedonline{24 February 2006}

\volumenumber{6}
\issuenumber{}
\publicationyear{2006}
\papernumber{3}
\startpage{91}
\endpage{117}

\doi{}
\MR{}
\Zbl{}

\makeatletter
\def\cnewtheorem#1[#2]#3{\newtheorem{#1}{#3}[section]
\expandafter\let\csname c@#1\endcsname\c@thm}

\makeatother  

\newtheorem{thm}{Theorem}[section]
\newtheorem{lem}[thm]{Lemma}
\newtheorem{cor}[thm]{Corollary}
\newtheorem{prop}[thm]{Proposition}

\DeclareMathOperator{\gq}{(RGQIE)} 
\DeclareMathOperator{\hd}{Hd}
\DeclareMathOperator{\nd}{Nbhd}
\DeclareMathOperator{\bold{SL_2}}{\bold{SL_2}}

\DeclareMathOperator{\rank}{rank} 

\DeclareMathOperator{\bold{SL_3}}{\bold{SL_3}}
\DeclareMathOperator{\bold{SL_n}}{\bold{SL_n}}

\newcommand{\ka}{\kappa}

\newcommand{\rt}{\rightarrow}

\newcommand{\ap}{{\mathcal{A}}}

\newcommand{\se}{\subseteq}
\newcommand{\e}{\varepsilon}

\newcommand{\ta}{K \gamma_{\mathfrak{p}}}
\newcommand{\de}{\delta}
\newcommand{\tgd}{K \gamma_{\mathfrak{p}}(\de)}
\newcommand{\p}{\phi}
\newcommand{\pn}{\phi ^{-1}}
\newcommand{\g}{\gamma}
\newcommand{\pgd}{\pi(\gr ,\de)}
\newcommand{\sa}{\measuredangle _}
\newcommand{\xy}{X_\infty}
\newcommand{\gy}{\gamma_\infty}
\newcommand{\py}{\pi_\infty}
\newcommand{\ey}{e_\infty}
\newcommand{\xr}{X_{\mathfrak{p}}}
\newcommand{\gr}{\gamma_{\mathfrak{p}}}
\newcommand{\pr}{\pi_{\mathfrak{p}}}
\newcommand{\er}{e_{\mathfrak{p}}}
\newcommand{\phr}{\phi_{\mathfrak{p}}}
\newcommand{\phy}{\phi_{\infty}}

\begin{document}

\begin{asciiabstract}
We give a new proof of a theorem of Kleiner-Leeb: that any
quasi-isometrically embedded Euclidean space in a product of symmetric
spaces and Euclidean buildings is contained in a metric neighborhood
of finitely many flats, as long as the rank of the Euclidean space is
not less than the rank of the target. A bound on the size of the
neighborhood and on the number of flats is determined by the size of
the quasi-isometry constants.  

Without using asymptotic cones, our proof focuses on the intrinsic
geometry of symmetric spaces and Euclidean buildings by extending the
proof of Eskin-Farb's quasiflat with holes theorem for symmetric
spaces with no Euclidean factors.
\end{asciiabstract}

\begin{abstract}
We give a new proof of a theorem of Kleiner--Leeb: that any
quasi-isometrically embedded Euclidean space in a product
of symmetric spaces and Euclidean buildings is contained in a
metric neighborhood of finitely many flats, as long as the rank of
the Euclidean space is not less than the rank of the target. A
bound on the size of the neighborhood and on the number of flats
is determined by the size of the quasi-isometry constants.

Without using asymptotic cones, our proof focuses on the intrinsic
geometry of symmetric spaces and Euclidean buildings by extending
the proof of Eskin--Farb's quasiflat with holes theorem for
symmetric spaces with no Euclidean factors.
\end{abstract}

\maketitle

\section{Introduction}\label{sec1}

We will give a new proof and a generalization of the following result:

\begin{thm}[Kleiner--Leeb]\label{1.1}
Let $\mathbb{E}^m$ be $m$--dimensional Euclidean space, and suppose
$\varphi \co \mathbb{E}^m \rt X$ is a $(\kappa , C)$
quasi-isometric embedding, where $X$ is a product of symmetric
spaces and Euclidean buildings and $m$ equals the rank of $X$.
Then there exist finitely many flats $F_1,F_2,\ldots,F_M \se X$ such
that
$$\varphi (\mathbb{E}^m)\se \nd _N \Big(\bigcup_{i=1}^M F_i\Big),$$
where $M=M(\kappa ,X)$ and $N=N(\kappa ,C ,X)$. \end{thm}

\fullref{1.1} was proved by Kleiner and Leeb in \cite{K-L}. It can be
used to give a new proof of a conjecture of Margulis from the 1970s
(also proved in \cite{K-L}) that any self-quasi-isometry of $X$ as
above is a bounded distance from an isometry when all factors
correspond to higher rank simple groups. For an indication as to how
\fullref{1.1} can be used to give a proof of this fact, see \cite{E-F1}
where Eskin--Farb give a proof of \fullref{1.1}  and Margulis' conjecture
in the case when $X$ is a symmetric space.

Our proof of \fullref{1.1} does not use asymptotic cones as the
proof of Kleiner--Leeb does. Rather, we adapt results of Eskin--Farb
who used large-scale homology to characterize quasiflats in
symmetric spaces without Euclidean factors in a way that allowed
for the absence of large regions in the domain of a quasiflat (a
``quasiflat with holes"). Thus, we provide a marriage between the
quasiflats theorems of Kleiner--Leeb and Eskin--Farb: a quasiflats
theorem that allows for products of symmetric spaces and Euclidean
buildings in the target of a quasiflat, and for holes in the
domain; see \fullref{1.2} below.  \fullref{1.1} occurs as a special
case.

Allowing for holes in our quasiflats leads to applications for the
study of the large-scale geometry of non-cocompact $S$--arithmetic
lattices; see Wortman \cite{W,T-W-W}.

\medskip
\textbf{Bibliographic note}\qua The full theorem of Kleiner--Leeb is more
general than \fullref{1.1} as it allows for generalized Euclidean
buildings in the target of $\varphi$. However, \fullref{1.1} does
include all of the standard Euclidean buildings that are naturally
acted on by reductive groups over local fields.

\medskip
\textbf{Quasiflats with holes}\qua For constants $\kappa \geq 1$ and
$C \geq 0$, a $(\kappa,C)$ \emph{quasi-isometric embedding} of a
metric space X into a metric space Y is a function $\varphi \co X
\rt Y$ such that for any $x_1 ,x_2 \in X$:
$$\frac{1}{\kappa}d(x_1,x_2)-C \leq d(f(x_1),f(x_2)) \leq \kappa d(x_1,x_2)+C.$$
For a subset of Euclidean space $\Omega \se \mathbb{E}^m$, we let
$$\Omega _{(\e , \rho)} = \{\, x\in \Omega \mid  B_y\big(\e d(x,y)\big)\cap \Omega \neq
\emptyset \mbox{ for all } y\in \mathbb{E}^m - B_x(\rho) \,\},$$
where we use the notation $B_z(r)$ to refer to the ball of radius
$r$ centered at $z$. Hence, $\Omega _{(\e , \rho)}$ is the set of
all points $x\in \Omega$ which can serve as an observation point
from which all points in $\mathbb{E}^m$ (that are a sufficient
distance from $x$) have a distance from $\Omega$ that is
proportional to their distance from $x$.

A special case to keep in mind is that if $\Omega = \mathbb{E}^m$,
then $\Omega_{(\e,\rho)} = \mathbb{E}^m$ for any $\e \geq 0$ and
$\rho \geq 0$.

A \emph{quasiflat with holes} is the image of $\Omega _{(\e ,
\rho)}$ under a quasi-isometric embedding $\phi \co\Omega \rt X$.

Before stating our main result, recall that for a metric space
$X$, the \emph{rank} of $X$ (or $\rank(X)$ for short) is the
maximal dimension of a flat in $X$. Now we have the following
generalization of \fullref{1.1}:

\begin{thm}[Quasiflats with holes]\label{1.2} Let $\varphi \co\Omega
\rt X$ be a $(\kappa , C)$ quasi-isometric embedding where $X$ is
a product of symmetric spaces and Euclidean buildings, $\Omega \se
\mathbb{E}^m$, and $m\geq \rank(X)$. There are constants
$M=M(\kappa ,X)$ and $\e _0 = \e _0 (\kappa ,X)$, such that if $\e
< \e_0$, then there exist flats $F_1,F_2,\ldots,F_M \se X$ such that
$$\varphi (\Omega _{(\e , \rho)})\se \nd _N \Big(\bigcup_{i=1}^M F_i\Big),$$
where $N=N(\kappa ,C , \rho ,X)$. \end{thm}

\textbf{Quasirank}\qua We remark that by comparing the volume of the
domain and image of a function $\varphi$ satisfying the hypotheses
of \fullref{1.2}, it is clear that no quasi-isometric embeddings
exist of a Euclidean space into $X$ when the dimension of the
Euclidean space is greater than the rank of $X$. This observation
is not new and follows very easily from the pre-existing
quasiflats theorems. However, we choose to state our theorem in
this more general manner since the proof given below does not
depend on the dimension of the Euclidean space once its dimension
at least equals the rank of $X$, and our proof will run more
smoothly if we allow for dimensions larger than the rank of $X$.

\medskip\textbf{Applications for quasiflats}\qua One would like to
characterize quasiflats as a starting point for understanding
quasi-isometries of a lattice as Mostow did for cocompact
lattices. 
(See Morse \cite{Mor}, Mostow \cite{Mo}, Pansu \cite{Pa},
Kleiner--Leeb\cite{K-L}, Eskin--Farb \cite{E-F1}, Eskin \cite{Es},
Wortman \cite{W,T-W-W} for the details of this brief sketch.)

The basic example of a quasiflats theorem is the Morse--Mostow
Lemma which states that a quasi-isometric embedding of
$\mathbb{R}$ into a rank one symmetric space has its image
contained in a metric neighborhood of a unique geodesic.

For general symmetric spaces and Euclidean buildings $X$, it is
not the case that a quasi-isometrically embedded Euclidean space
is necessarily  contained in the neighborhood of a single flat.
(Recall that a \emph{flat} is an isometrically embedded Euclidean
space.) If, however, the dimension of a quasi-isometrically
embedded Euclidean space is equal to the dimension of a maximal
flat in $X$, then its image will be contained in a neighborhood of
finitely many flats.

Quasiflats can be used in the study of quasi-isometries of
cocompact lattices as follows. First, we may assume that any
self-quasi-isometry of a cocompact lattice in a semisimple Lie
group is a quasi-isometry of its orbit in an appropriate product
of symmetric spaces and Euclidean buildings, $X$. Second, since
any flat in $X$ is necessarily contained in a metric neighborhood
of the cocompact lattice orbit, we can restrict the quasi-isometry
to any flat and examine its image. The space $X$ has a boundary at
infinity which is defined in terms of the asymptotic behavior of
flats, so in determining the images of flats we are finding a map
on the boundary of $X$. Finally---as long as $X$ contains no
factors that are real hyperbolic spaces, complex hyperbolic
spaces, or trees---one can deduce from the properties of the
boundary map that the quasi-isometry is a finite distance from an
isometry.

The story is different for non-cocompact lattices. Generic flats
in $X$ will not be contained in a neighborhood of a non-cocompact
lattice orbit. Hence, we cannot apply the same proof technique.

However, the generic flat will have a substantial portion of its
volume contained in a neighborhood of a non-cocompact lattice
orbit. With an eye towards this feature, Eskin--Farb provided a
foundational tool for studying quasi-isometries of non-cocompact
lattices in real semisimple Lie groups by defining and
characterizing quasiflats with holes in symmetric spaces.

Using quasiflats with holes in symmetric spaces, Eskin developed a
boundary map in the non-cocompact lattice case for real groups en
route to proving that any quasi-isometry of a higher rank
arithmetic group is a finite distance from a commensurator.

By allowing for Euclidean building factors in the image of a
quasiflat with holes, we will be able to use this same approach to
analyze quasi-isometries of non-cocompact lattices in semisimple
Lie groups over arbitrary local fields.

\medskip\textbf{Outline}\qua Our proof of \fullref{1.2} in the case that $X$ is
a Euclidean building is self-contained aside from results of
Eskin--Farb on the large-scale homology of pinched sets in
Euclidean space and some consequences of those results. Hopefully,
the reader who is interested in only the case when $X$ is a
building can read through our proof without having to consider
symmetric spaces.

In the general case, when $X$ is a nontrivial product of a
symmetric space and a Euclidean building, we rely heavily on the
results of Eskin--Farb for symmetric spaces. Our approach is to
project the quasiflat with holes into the building factor $\xr$,
and into the symmetric space factor $\xy$. By projecting the
quasiflat with holes to $\xr$, we can apply arguments below that
were created expressly for buildings while ignoring the symmetric
space factor. Conversely, by projecting the quasiflat with holes
to $\xy$, we can directly apply most of the content of \cite{E-F1} 
to analyze the image. After examining the image in each factor,
we piece together the information obtained in the full space $X$
to obtain our result.

Thus, in our approach to proving \fullref{1.2}, we will try to avoid
dealing with the product space $X$. We do this since arguments for
symmetric spaces and Euclidean buildings (although extremely
similar in spirit) have to be dealt with using different tools.

The approach of projecting to factors is taken from the work of
Eskin--Farb as well. Their test case for their general theorem was
when $X=\mathbb{H}^2 \times \mathbb{H}^2$, and they used the
projection method to reduce most of the proof to arguments in the
hyperbolic plane \cite{E-F2}.

In \fullref{sec2} we will show that certain subspaces in $X$ which
behave like rank one spaces cannot accommodate quasi-isometric
embeddings of large Euclidean sets. This fact will be formulated
more precisely in terms of homology.

Some of the nearly rank one spaces are then glued together to give
a ``degenerate space" in $X$ which is a fattening of the singular
directions in $X$ with respect to a given basepoint. (Recall that
a direction is singular if it is contained in more than one flat.)
Using a Mayer--Vietoris sequence, it can be shown that the
degenerate space cannot accommodate quasi-isometric embeddings of
large Euclidean sets of large dimension. It is at this point where
we apply our hypothesis that the dimension of $\mathbb{E}^m$
equals, or exceeds, the rank of $X$.

In \fullref{sec3} we begin to analyze the asymptotic behavior of
quasiflats with holes. We define---following Eskin--Farb---what it
means for a direction in a quasiflat with holes to limit on a
point in the boundary at infinity of $X$.

The results of \fullref{sec2} show that the image of a quasiflat with
holes must have a substantial intersection with the complement of
the degenerate space. (The complement of the degenerate space is
the region of $X$ for which limit points are defined.) We argue
further to show that limit points exist.

Since the nondegenerate space behaves much like a rank one space
itself, we can show that the image of a quasiflat with holes in
the nondegenerate space cannot extend in too many directions (i.e.\
the number of limit points is bounded). We construct our bound by
contrasting the polynomial growth of Euclidean space with the high
cost of travelling out in different directions in a rank one
space. It is from the finite set of limit points that the finite
set of flats from the conclusion of \fullref{1.2} is constructed.

\fullref{sec4} contains a few lemmas to insure that all definitions
depending on basepoints are well-defined up to a constant.

We conclude in \fullref{sec5} with a proof of \fullref{1.2}. Results from
Sections \ref{sec2}, \ref{sec3} and \ref{sec4} are used in the proof.

\medskip\textbf{Definitions}\qua Recall that a \emph{polysimplex}
 is a product of simplices. Replacing simplices with polysimplices
 in the definition of a simplicial complex creates what is called
 a \emph{polysimplicial complex}.

 A \emph{Euclidean building} $\xr$ is a polysimplicial
 complex endowed with a metric $d_{\mathfrak{p}}$ that satisfies
  the four properties below:

\begin{itemize}
\item[(i)]There is a family, $\{ \mathcal{A}_\alpha \}$, of
subcomplexes of $\xr$ such that each $\ap _\alpha$ is isometric to
$\mathbb{E}^{\text{dim}(\xr)}$ and $\xr = \bigcup _\alpha
\mathcal{A}_\alpha$. Each $\mathcal{A}_\alpha$ is called an
\emph{apartment}.

\item[(ii)]Any two polysimplices of maximal dimension (called
\emph{chambers}) are contained in some $\mathcal{A}_\alpha$.

\item[(ii)]If $\mathcal{A}_\alpha$ and $\mathcal{A}_\beta$
are two apartments each containing the chambers $\mathfrak{c}_1$
and $\mathfrak{c}_2$, then there is an isometric polysimplicial
automorphism of $X$ sending $\mathcal{A}_\alpha$ to
$\mathcal{A}_\beta$, and fixing $\mathfrak{c}_1$ and
$\mathfrak{c}_2$ pointwise.

\item[(iv)]The group of isometric polysimplicial automorphisms
of $\xr$ acts transitively on the set of chambers.
\end{itemize}

Note that condition (iv) is nonstandard. Often one assumes the
stronger condition that a building be \emph{thick}. We desire to
weaken the thickness condition to condition (iv) so that Euclidean
space can naturally be given the structure of a Euclidean
building.

Also notice that we do not assume $X_{\mathfrak{p}}$ to be locally
finite. Hence, we are including the buildings for, say,
$\mathbf{GL_n} (\mathbb{C}(t))$ in our examination.

Along with the nonstandard definition of a Euclidean building
given above, we also give the standard definition of a
\emph{symmetric space} as a Riemannian manifold $\xy$ such that
for every $p \in \xy$, there is an isometry $g$ of $\xy$ such that
$g(p)=p$ and the derivative of $g$ at $p$ equals $-\text{Id}$.

\medskip\textbf{Conventions}\qua Throughout this paper we will be examining
products of symmetric spaces and Euclidean buildings. Since
Euclidean space is a Euclidean building by our definition, we may
assume that our symmetric spaces do not have Euclidean factors.
This will allow us to more readily apply results from \cite{E-F1}
where it is assumed that the symmetric spaces have no Euclidean
factors.

We may also assume that our symmetric spaces do not have compact
factors. Otherwise we could simply compose the quasi-isometry
$\varphi$ from \fullref{1.2} with a projection map to eliminate the
compact factors, then apply \fullref{1.2}, pull back the flats
obtained to the entire symmetric space, and increase the size of
$N$ by the diameter of the compact factors.

\medskip\textbf{Notation}\qua If $a$ and $b$ are positive numbers we write
$a\ll b$ when there is a constant $\lambda=\lambda(X, \kappa)<1$
such that $a<\lambda b$. If there are variables $x_1,\ldots,x_n$ and a
constant $\eta=\eta (X,\kappa, x_1,\ldots,x_n)<1$ such that $a<\eta
b$, then we write $a \ll _{(x_1,\ldots,x_n)} b$. We will use the
notation $a=O(b)$ to mean that $a<\lambda b$ for some constant
$\lambda=\lambda(X, \kappa)$ without specifying the size of
$\lambda$.

\medskip\textbf{Remarks}\qua With modification to only the conclusion of the
proof of \fullref{3.6}, our results hold when $\mathbb{E}^m$ is
replaced by a $1$--connected nilpotent real Lie group. For example,
this shows that a Heisenberg group cannot quasi-isometrically
embed into $\mathbf{SL_4}(k)$ for any locally compact nondiscrete
field $k$.

Also the proof presented below can be modified in \fullref{3.2} to
allow for the presence of $\mathbb{R}$--buildings in the target of
the quasiflat with holes.

\medskip\textbf{Acknowledgements}\qua Benson Farb was my PhD thesis advisor
under whose direction this work was carried out. I thank him for
suggesting this problem to me, and for his constant support and
encouragement.
Thanks also to Alex Eskin for listening to many of my ideas and
for providing feedback.
Thanks to Tara Brendle, Dan Margalit, Karen Vogtmann and a
referee for valuable comments made on an earlier draft.
I would also like to thank the University of Chicago for
supporting me as a graduate student while I developed the ideas in
this paper, and Cornell University for the pleasant working
environment given to me while I completed the writing of this
paper.
I was supported in part by an NSF Postdoctoral Fellowship.

\section{Pinching functions and homology}\label{sec2}

Throughout the remainder, let $\xr$ be a Euclidean building with a
chosen basepoint $\er \in \xr$, and let $\xy$ be a symmetric space
with basepoint $\ey \in \xy$. We will assume that $\xy$ has no
compact or Euclidean factors (see the conventions in the preceding
section).

We let $X =\xy \times \xr$, and we define $\pi _\infty \co X \rt
\xy$ and $\pi _{\mathfrak{p}}\co X \rt \xr$ to be the projection
maps. Define the point $e \in X$ as the pair $(\ey, \er)$.

Throughout we let $n \in \mathbb{N}$ equal $\rank(X)$.

\medskip\textbf{Graded quasi-isometric embeddings}\qua We will put quasiflats
with holes aside until the final section of this paper. We
concentrate instead on embeddings of entire Euclidean spaces into
$X$ under a weaker assumption than our map is a quasi-isometry.

For points $x,y_1,y_2,\ldots,y_n \in X$ and a number $\rho \geq 0$,
we let
$$D_x ( \rho ; y_1, y_2, \ldots ,y_n ) = \max \{ \rho, d(x, y_1),
\ldots, d(x,y_n)\}.$$ For numbers $\kappa \geq 1$, $\rho \geq 0$, and
$\varepsilon \geq 0$, we define a function $\p \co X \rightarrow
Y$ to be a $(\ka ,\rho,\e)$ \emph{graded quasi-isometric embedding
based at} $x\in X$ if for all $z,w \in X$:
\[ \frac{1}{\ka} d(z,w)-\e D_x(\rho;z,w) \leq d(\p (z),
\p(w) ) \leq \ka d(z,w)+\e D_x(\rho;z,w).
\]
A function $\p \co X \rightarrow Y$ is called $(\kappa ,\rho)$
\emph{radial} at $x\in X$ if for all $z\in X$:
\[ \frac{1}{2 \kappa} D_x(\rho;z) \leq d(\p (z),
\p(x) ) \leq (2 \kappa ) D_x(\rho;z).
\]
Combining the two definitions above, $\p \co  X \rt Y$ is a $
(\kappa, \rho, \e )$ \emph{radial graded quasi-isometric
embedding} ($\gq$ for short) based at $x$ if it is a $(\kappa,
\rho, \e )$ graded quasi-isometric embedding at $x$, and $\kappa$
radial at $x$.

In the proof of \fullref{1.2}, we will see that one can easily
extend the domain of a quasiflat with holes to all of
$\mathbb{E}^m$ in such a way that the extension is a $\gq$. From
the behavior of $\gq$'s that is characterized in Sections \ref{sec2}
through \ref{sec4}, we will be able to characterize the image of a
quasiflat with holes.

Until explicitly stated otherwise, let $\p \co  \mathbb{E}^m \rt
X$ be a $(\kappa, \e , \rho)$ $\gq$ based at $0$ with $\p (0) =e$.
The image of such a function is a \emph{graded quasiflat}.

\medskip\textbf{Pinching on rays in buildings}\qua Let
$$K=\{\, g \in \text{Isom}(X) \mid ge=e \,\},$$ and let $\gr
\co [0,\infty ) \rt \{\ey\}\times \xr$ be a geodesic ray with $\g
_{\mathfrak{p}} (0)=e$. The space $\ta$ is a topological tree as can
be seen by restricting the geodesic retraction $\xr \rt \{\er\}$.
However, the tree $\ta$ will often not be convex. These trees in
$X$ are negatively curved, and our first goal is to show that
large subsets of Euclidean space cannot embed into them, or even
into small enough neighborhoods of them. This in itself is
straightforward to show, but we shall want to handle this problem
in a way that allows us to conclude that large Euclidean sets
cannot embed into fattened neighborhoods of $K$ translates of
certain $(n-1)$--dimensional spaces.

 Let $$\tgd=\{\, x\in \{\ey \} \times \xr \mid d(x,t) < \de d(x,e) \mbox{
for some } t\in \ta \,\},$$ so that $\tgd$ is a neighborhood of
$\ta$ in $\{\ey\}\times \xr $ that is fattened in proportion to
the distance from the origin by a factor of $\delta$. We will want
to project $\tgd$ onto $\ta$ where calculations can be made more
easily.

Define $$\pgd \co  \tgd \rt \ta$$ by choosing for any $x\in \tgd$,
some $\pgd (x) \in \ta$, such that $$d(x,\pgd (x))\leq \de
d(x,e).$$ By definition, $\pi (\gr , \de)$ only modifies distances
by a linear error of $\de$, so composing with $\p$ will still be a
$\gq$. Precisely, we have the following:

\begin{lem}\label{2.1} 
If $\e < \de < 1/2$, then $\pgd \circ \p \co  \pn (\tgd) \rt
\ta$ is a $(2 \ka,\rho,5\ka \de)$ $\gq$ based at $0$. \end{lem}

\begin{proof}
 Verifying that $\pgd \circ \p$ is a graded
quasi-isometric embedding is an easy sequence of inequalities:
\begin{align*}
d\big(\pi(\gr& ,\de)\, \circ\,   \p (x)\,,\, \pgd \circ \p (y)\big) \\
& \leq d\big(\pgd \circ \p (x)\,,\,
 \p (x)\big)+d\big(\pgd \circ \p (y)\,,\, \p (y)\big) + d\big(\p (x)\,,\, \p (y)\big) \\
 & \leq d\big(\p(x)\,,\, \p(y)\big)+2\de D_e\big(0;\p (x),\p(y)\big) \\
 & \leq \ka  d\big(x,y\big) + \e D_0 \big(\rho ;x,y\big) + 4\ka
     \de D_0\big(\rho;x,y\big).
\end{align*}
The other inequality is similar.

That $\pgd \circ \p$ is radial is also straightforward:
\begin{align*} d\big(\pgd \circ \p (x)\,,\, e\big) & \leq  d\big(\pgd \circ \p
(x)\,,\, \p (x)\big)+d\big(\p(x)\,,\,e\big) \\
& \leq (1+ \de )d\big(\p (x)\,,\,e\big) \\
& \leq 2\ka (1+ \de) D_0\big(\rho ;x\big).
\end{align*}
Again, the other inequality is similar.
\end{proof}

As in \cite{E-F1}, for numbers $r \geq 0$, $\eta >1$, and $\beta
>0$, we define an $(r,\eta,\beta)$ \emph{pinching function} on a set
$W\subseteq \mathbb{E}^m$  to be a proper, continuous function
$f\co W \rightarrow \mathbb{R}_{\geq 0}$ such that for any $x,y
\in W$, we have $d(x,y)<\beta s$ whenever the following two
properties hold:
\begin{itemize}
\item[(i)]$r \leq s \leq f(x) \leq f(y) \leq \eta s$;

\item[(ii)] there is a path $\psi \co [0,1] \rightarrow
W$ such that $\psi (0)=x$, $\psi(1)=y$, and
$s\leq f(\psi (t))$ for all $t\in [0,1]$.
\end{itemize}

If there exists an $(r,\eta,\beta)$ pinching function on some $W
\subseteq \mathbb{E}^m$, then we say that $W$ is
$(r,\eta,\beta)$--\emph{pinched}.

Eskin--Farb used pinching functions as a means of showing that
large Euclidean sets cannot quasi-isometrically embed into certain
negatively curved subspaces of symmetric spaces. To show the
analogous result for our general $X$, we will first construct a
pinching function for $\pn (\tgd )$. Since Eskin--Farb constructed
a pinching function on the similarly defined sets $\pn (K \g
_\infty (\de))$, we will then be in a position to handle the case
for a general ray by pulling back pinching functions obtained
through projection to factors.

Our candidate for a pinching function on $\pn (\tgd )$ is $$f(\gr
, \de)\co  \pn (\tgd) \rt \mathbb{R}_{\geq 0},$$ $$f(\gr ,
\de)(x)=d(\pgd \circ \p (x), e).\leqno{\hbox{where}} $$

\begin{lem}\label{2.2} 
If $\e < \de <1/2$, then $f(\gr , \de)$ is a $(5 \ka \rho,
1+\de , 84 \ka ^3 \de )$ pinching function on the set $\pn (\tgd)
\se \mathbb{E}^m$. \end{lem}

 \begin{proof} Note that we may assume $\pgd \circ \p$ is continuous by
  a connect-the-dots argument. Hence, $f(\gr , \de)$ is clearly continuous and
proper. We assume $x,y \in  \pn (\tgd)$ are such that $$5 \ka \rho
\leq s \leq f(\gr , \de)(x) \leq f(\gr , \de)(y) \leq (1+\de )s,$$
and there is a path $\psi \co [0,1] \rt \pn (K\gr (\de ))$ with $s
\leq f(\gr , \de)(\psi (t))$ for all $t \in [0,1]$.

By the radial condition of \fullref{2.1}, $$5 \ka \rho \leq d\big(\pgd
\circ \p (x)\,,\,e\big)\leq 4\ka D_0(\rho ;x).$$ It follows that
$\rho < d(x,0)$. Hence, by the radial condition of \fullref{2.1} and
our pinching assumptions,
$$d(x,0) \leq 4 \ka d\big(\pgd \circ \p (x)\,,\,e\big) \leq
4 \ka (1+ \de )s.$$
 The existence of $\psi$ implies that $\pgd \circ \p
(x)$ and  $\pgd \circ \p (y)$ are in the same connected component
of $\ta - B_e (s)$. Therefore, $$d\big( \pgd \circ \p (x)\,,\,
\pgd \circ \p (y)\big) \leq 2\de s.$$ We may assume $d(x,0)\geq
d(y,0)$. Then, by the graded condition of \fullref{2.1},
$$2\de s \geq \frac{1}{2 \ka}d(x,y) - (5\ka \de )d(x,0) \geq
\frac{1}{2 \ka }d(x,y)-(5 \kappa \de ) 4 \kappa (1+ \de ) s.$$
That is, $d(x,y) < 84 \ka ^3 \de s$.
\end{proof}

\textbf{Graded neighborhoods}\qua For a set $Y \se X$,
 we can create a neighborhood of $Y$ by fattening points in $Y$ in
 $\de$--proportion to their distance from $e$. In symbols, we let
$$Y [\de ]=\{\, x\in X \mid d(x,y) < \de d(x,e) \mbox{ for some }
y\in Y \,\}.$$

\medskip\textbf{Pinching on general rays}\qua Lemma 6.8 in \cite{E-F2}
demonstrates a pinching function for sets of the form $\pn (K \gy
(\de ))$ where $\gy \co  [0, \infty ) \rightarrow \xy \times \{\er
\}$ is a geodesic ray, and $K \gy (\de ) \se \xy \times \{ \er \}$
is defined analogously to $K \gr (\de ) \se \{ \ey \} \times \xr$.
We can use this pinching function along with the pinching function
from \fullref{2.2} to show that $\pn (K \g [\de ])$ is a pinched set,
where $\g \co [0, \infty ) \rightarrow X$ is an arbitrary geodesic
ray with $\g (0)=e$. Our argument proceeds by simply applying our
already existing pinching functions to the image of $K\g [\de ]$
under the projection maps onto the factors of $X$.

We want to define a real valued \emph{tilt parameter}, $\tau$, on
the space of geodesic rays $\g \co  [0 ,\infty ) \rightarrow X$
with $\g(0)=e$. The parameter will measure whether $\g$ leans more
towards the $\xr$ or the $\xy$ factor. Notice that any such $\g$
can be decomposed as $\g (t) = (\gy (t), \gr(at))$ for some number
$a\geq 0$, and all $t\geq 0$, where $\gy \subseteq \xy$ and $\gr
\subseteq \xr$ are unit speed geodesic rays based at $\ey$ and
$\er$ respectively. Now we simply set $\tau (\g) =a$. (For $\tau$
to be defined everywhere we
 allow for the case when $a=\infty$, which is just to say that
$\g$ is contained in the building factor.) Hence, if $\tau (\g)
> 1$ (resp.\ $< 1$) then $\g$ is leaning towards the building
factor (resp.\ symmetric space factor), and when creating a
pinching function on $K \g [\de]$ it will be most efficient to
project onto the $\xr$ (resp.\ $\xy$) factor of $X$.

We begin with the following technical observation.

\begin{lem}\label{2.3} 
Assume $\g \co  [0 , \infty ) \rightarrow X$
 is a geodesic ray with
$\g(0)=e$ and that $y\in K\g [\de]$. Then,

\begin{itemize}
\item[(i)] $\pr(y)\in K\gr \Big( \delta
\sqrt{1+\cot^2(|\tan^{-1}\tau(\g)-\sin^{-1}\delta|_{+})}\; \Big)$,
 \quad and

\item[(ii)]  $\py(y)\in K\gy \Big(\delta
\sqrt{1+\cot^2(|\tan^{-1}1/\tau(\g)-\sin^{-1}\delta|_{+})}\;
\Big)$,

\end{itemize}

where $|x| _+ = \max \{ x,0 \}$. \end{lem}

 \begin{proof} By definition of $K \g [\delta]$
there exists a $t \geq 0$ and a $k\in K$ such that
\begin{align*}
d\big(\pr(y)\,,\,k\gr(\tau(\gamma)t)\big) & = d\big(\pi
_{\mathfrak{p}}(y)\,,\,
\pi_{\mathfrak{p}}(k\gamma (t))\big) \\
& \leq d\big(y\,,\,k\g(t)\big) \\
 & <\delta d(y,e) \\
 & \leq \delta \sqrt{d(\pr (y),\er)^2 + d(\py(y),\ey)^2}.
 \end{align*}
Using straightforward trigonometry it can be verified that
$$d(\py(y),\ey)\leq d(\pr(y),\er)
\cot(|\tan^{-1}\tau(\g)-\sin^{-1}\delta|_{+}).$$ Then (i) follows.
The proof of (ii) is similar.
\end{proof}

We will use part (i) of the previous lemma to create a pinching
function for geodesic rays that tilt towards $\xr$. This is the
content of \fullref{2.5}, but we will first note that the projection
onto $\xr$ does not significantly distort distances.

\begin{lem}\label{2.4} Let $\g\co [0,\infty) \rt X$ be
 a geodesic ray with $\g(0)=e$. If $\e<\de$ and $\tau (\g )\geq 1$, then $\pr
  \circ \p \co  \pn (K\g[\de ])\rightarrow \xr$ is
   a $(2\kappa, \rho, \eta _1)$ $\gq$ where $\eta _1=O(\de )$.
   \end{lem}

\begin{proof} Note that on $ K\g[\de ]$, $\pr$ is a  $(2,0,O(\de ))$ $\gq$
where $2$ is an upper bound given by our restriction on $\tau (\g
)$. Composition with $\p$ completes the result.
\end{proof}

Now for the pinching function:

\begin{lem}\label{2.5} Let $\g\co  [0,\infty) \rt X$ be
 a geodesic ray with $\g(0)=e$. For $\tau (\g )\geq 1$ and $\e<\de \ll 1$, the set
  $\pn (K\g [\de ]) \se \mathbb{E}^m$ is
$(10 \kappa \rho, 1+\de , O(\de ) )$--pinched. \end{lem}

\begin{proof} Let $\de _{\mathfrak{p}} =\max \Big \{ 2 \eta _1 ,\de
\sqrt{1+\cot^2(\tan^{-1}\tau (\g )-\sin^{-1}\delta)} \Big \} $,
and note that our conditions on $\tau (\g )$ and $\de$ imply that,
say,
$$1<\sqrt{1+\cot^2(\tan^{-1}\tau (\g )-\sin^{-1}\delta)}<2.$$
By \fullref{2.3}, $\pr (K\g [\de ])\subseteq K \gr (\de
_{\mathfrak{p}})$. Hence, we can choose our pinching function
$g\co \pn (K\g [\de ]) \rightarrow \mathbb{R}_{\geq 0}$ to be
given by
$$g(z)=d(\pi (\gr , \de _{\mathfrak{p}})\circ \pr \circ \p
(z),\er).$$ Indeed, we can use \fullref{2.4} to replace $\p$ with $\pr
\circ \p$ in \fullref{2.2}. It follows that $g$ is a $(10\kappa \rho,
1+\de _{\mathfrak{p}}, 672 \kappa ^3 \de _{\mathfrak{p}})$
pinching function.
\end{proof}

If $\tau (\g) \leq 1$, we can apply Lemma 2.3 to Lemma 6.8 of
\cite{E-F1} and obtain a similar result. Hence, we have a
pinching function on $\pn (K \g [\de ])$ for any geodesic ray $\g$
that is based at the origin. Precisely, we have the following:

\begin{lem}\label{2.6} If $\e \ll \de \ll 1$, then the set $\pn (K \g [\de ])\se
\mathbb{E}^m$ is $(r_0,1+O(\de),O(\de))$--pinched for any geodesic
ray $\g\co [0,\infty) \rt X$ with $\g(0)=e$. Here
$r_0=r_0(X,\kappa, \rho, \de)$. \end{lem}

\textbf{Homology results of Eskin--Farb and their consequences}\qua
Pinching functions were introduced in \cite{E-F1} as a tool for
showing that sets which simultaneously support Euclidean metrics
and ``quasinegatively curved" metrics must be small and, hence,
cannot have any interesting large-scale homology. Precisely, we
can use our \fullref{2.6} in the proof of Corollary 6.9 from \cite{E-F1} 
to show:

\begin{lem}\label{2.7} There exists a $\nu_1>0$ such that if $1 \ll_{(\rho, \de,
\e)} r$, while $\e\ll\de\ll1$ and $W\se \pn (K\g[\de])$, then the
homology of the inclusion map $\iota_* \co H_p(W \cup B_0(r)) \rt
H_p (W[\nu_1 \de ] \cup B_0(r))$ is zero for all $p\geq 1$.
\end{lem}

The above lemma can be used to show, for example, that the image
of $\p$ cannot be contained in $K\g[\de]$. Otherwise we could take
a sphere of large radius in place of $W$ to arrive at a
contradiction. This is an interesting fact, but we care to know
more. We are able to use this lemma to tell us that there are much
larger subspaces of $X$ that spheres cannot embed into.

The larger subspaces are defined in terms of walls, so we begin by
defining the latter. A subset $H\se X$ is called a \emph{wall} if
it is a codimension $1$ affine subspace of a flat that is
contained in at least two distinct flats. Note that the walls
through the point $e\in X$ comprise the singular directions from
$e$.

Our space $X$ resembles a rank one space, from the vantage point
of $e\in X$, in the regions bounded away from the singular
directions. Properties of negative curvature are a powerful tool,
so we will want to show the image of $\p$ has a substantial
portion of its image bounded away from the singular directions.

It is time to define $X_e(\de)$ as the $\de$--\emph{nondegenerate
space} at $e\in X$ consisting of those points in $X$ that are not
contained in any $\de$--graded neighborhood of a wall containing
$e$. That is $$X_e(\de)=\bigcap_{H \in \mathcal{W}_e}(H[\de])^c,$$
where $\mathcal{W}_e$ is the set of walls in $X$ that contain $e$.

The complement $X_e(\de)^c$ of the $\de$--nondegenerate space is
the $\de$--\emph{degenerate space}. We could repeat the definition
for the special case that $X$ is either a Euclidean building or a
symmetric space and obtain the sets $X_{\mathfrak{p},\er}(\de)$,
$X_{\mathfrak{p},\er}(\de)^c$, $X_{\infty,\ey}(\de)$, and
$X_{\infty,\ey}(\de)^c$.

Our goal for this section is to show that the image of $\p$ is
forced to travel in $X_e(\de)$. We can use \fullref{2.7} along with a
Mayer--Vietoris sequence to show that the image under $\p$ of very
large subsets of $\mathbb{E}^m$ indeed cannot be contained in $X_e
(\de )^c$. Note that in the Tits boundary of $X$, $X_e (\de )^c$
appears as a neighborhood of the $(n-2)$--skeleton. The spaces of
the form $K\g[\de]$ that we considered previously appear as
neighborhoods of a family of points in the Tits building. It is
clear how one would want to use \fullref{2.7} and a Mayer--Vietoris
argument to arrive at the following:

\begin{lem}\label{2.8} There exists a constant $\nu _2>0$, such that if
$1\ll_{(\rho, \de, \e)} r$ while $\e\ll\de\ll1$ and $W\se \pn (X_e
(\de )^c)$, then the homology of the inclusion map $\iota_* \co
H_{p}(W \cup B_0(r)) \rt H_{p} (W[\nu_2 \de ] \cup B_0(r))$ is
zero for all $p \geq n-1$. \end{lem}

The basic idea of the proof is clear but there are some
technicalities to consider. This is essentially Lemma 5.6 of
\cite{E-F1}, whose proof takes place in the Tits boundary where
there is no difference between symmetric spaces and buildings.
Hence, the proof carries over completely to prove our \fullref{2.8}.

\medskip\textbf{Unbounded, nondegenerate components of graded quasiflats}\qua
Note that the above lemma tells us that large metric
$(n-1)$--spheres in $\mathbb{E}^m$ cannot map into $X_e (\de )^c$
under $\p$. In Lemma 5.8 of \cite{E-F1}, this idea is extended to
show that unbounded portions of $\mathbb{E}^m$ map into $X_e (\de
)$ under $\p$. The arguments there only involve an application of
what is our \fullref{2.8} to the homology of Euclidean sets. The proof
applies verbatim to yield:

\begin{cor}\label{2.9} There is a constant
 $\nu _3>1$,
such that if $\e \ll \de \ll 1$ and $z \in \pn (X_e(\de ))$ with
$1 \ll _{(\de , \e ,\rho)} r \leq d(z,0)$, then the connected
component of $\pn (X_e (\de/ \nu_3))\cap B_0 (r)^c$ that contains
$z$ is unbounded. \end{cor}

\fullref{2.8} and \fullref{2.9} are the only results from this section
that will be used in the remainder of this paper. We will apply
\fullref{2.8} in \fullref{sec5} during the proof of \fullref{1.2}. 
\fullref{2.9} is used in the proof of \fullref{3.5} below to create a path in
the graded quasiflat that avoids the nondegenerate space and
accumulates on a point in the boundary of $X$.

\section{Limit points in Euclidean buildings}\label{sec3}

\textbf{Boundary metric}\qua A subset of a Euclidean building
$\mathfrak{S}\se \xr$ is called a \emph{sector} based at $x \in
\xr$, if it is the closure of a connected component of an
apartment less all the walls containing $x$.

Let $\what{X}_{\mathfrak{p}}$ be the set of all sectors based at
$\er$. For any $\mathfrak{S} \in \what{X}_{\mathfrak{p}}$, let
$\g _{\mathfrak{S}} \co [0, \infty ) \rt \mathfrak{S}$ be the
geodesic ray such that $\g _{\mathfrak{S}}(0) =\er$, and such that
$\g _{\mathfrak{S}} (\infty )$ is the center of mass of the boundary
at infinity of $\mathfrak{S}$ with its usual spherical metric. We
will also use $\gamma _{\mathfrak{S}}$ to denote the image of
$\gamma_{\mathfrak{S}} \co [0,\infty) \rt \mathfrak{S}$.

We endow $\what{X}_{\mathfrak{p}}$ with the metric
$\what{d}_{\mathfrak{p}}$ where
\begin{equation*}
\what{d}_{\mathfrak{p}} (\mathfrak{Y},\mathfrak{Z}) =
\begin{cases}
\, \quad \pi, &\text{if $\g _{\mathfrak{Y}} \cap \g _{\mathfrak{Z}} = \{ \er \}$;}\\
\frac{1}{|\g _{\mathfrak{Y}} \cap \g _{\mathfrak{Z}} |},
&\text{otherwise.}
\end{cases}
\end{equation*}
In the above, $|\g _{\mathfrak{Y}} \cap \g _{\mathfrak{Z}} |$ is the
length of the geodesic segment $\g _{\mathfrak{Y}} \cap \g
_{\mathfrak{Z}} $.

Note that $\what{d}_{\mathfrak{p}}$ is invariant under the action
of the stabilizer of $\er$ and is a complete ultrametric on
$\what{X}_{\mathfrak{p}}$. That $\what{d} _{\mathfrak{p}}$ is an
ultrametric means that it is a metric, and
$$\what{d}_{\mathfrak{p}} (\mathfrak{Y},\mathfrak{Z}) \leq \max
\{ \what{d}_{\mathfrak{p}} (\mathfrak{Y},\mathfrak{X}) ,
\what{d}_{\mathfrak{p}}(\mathfrak{X},\mathfrak{Z}) \} \text{
\quad for any }\mathfrak{Y,Z,X} \in \what{X}_{\mathfrak{p}}.$$ We
will use at times  that $$\mathfrak{Z}\in B_{\mathfrak{S}}(r) \text{
implies } B_{\mathfrak{Z}}(r) =B_{\mathfrak{S}}(r),$$ which is a
reformulation of the ultrametric property.

\medskip\textbf{Measuring angles}\qua We also introduce a notion of angle
between two points in a building as measured from $\er$. We first
define $\Phi _{\mathfrak{p}} \co \xr \rt
\mathcal{P}(\what{X}_{\mathfrak{p}})$ by $$\Phi _{\mathfrak{p}}
(x)= \{\, \mathfrak{S} \in \what{X}_{\mathfrak{p}} \mid x\in
\mathfrak{S} \,\},$$ where $\mathcal{P}(\what{X}_{\mathfrak{p}})$
denotes the power set of $\what{X}_{\mathfrak{p}}$.

Then for any $x,y \in \xr,$ we define
$$ \Theta _{\mathfrak{p}} (x ,y)=\inf \big\{ \, \what{d} _{\mathfrak{p}} (\mathfrak{S}_x,
\mathfrak{S}_y) \mid \mathfrak{S}_x \in \Phi _{\mathfrak{p}} (x)
\mbox{ and } \mathfrak{S}_y \in \Phi _{\mathfrak{p}} (y)
 \, \big\}.$$ We think of $ \Theta _{\mathfrak{p}} (x ,y)$ as measuring an angle between $x$
and $y$.

We will also be measuring angles formed by triangles in a single
apartment. Since apartments are Euclidean spaces, we can simply
use the Euclidean measure of angle. If $\ap \se \xr$ is an
apartment and $x,y,z\in \ap$, we let $\measuredangle _z ^\ap
(x,y)$ be the standard Euclidean angle in $\ap$ between $x$ and
$y$ as measured at $z$. For any subset $H \se \ap$, and points
$x,z \in \ap$, we let $$\measuredangle _z ^\ap (x,H)=\min
\{\measuredangle _z ^\ap (x,h) | h\in H \}.$$

\medskip\textbf{Core of a sector}\qua From here on we will assume that $0
\leq \de \leq 1$. For any $\mathfrak{S} \in
\what{X}_{\mathfrak{p}}$, we let $$\mathfrak{S}(\de )=\{\, x \in
\mathfrak{S} \mid d(\partial \mathfrak{S}, x)\geq \de d(e ,x)
\,\}.$$ We refer to $\mathfrak{S}(\de)$ as the $\de$--\emph{core}
of $\mathfrak{S}$. Note that
$$\bigcup _{\mathfrak{S} \in \what{X}_{\mathfrak{p}}}
\mathfrak{S}(\de ) = X_{\mathfrak{p},\er}(\de),$$ where
$X_{\mathfrak{p},\er}(\de)$ is the $\de$--nondegenerate space of
$\xr$ at $\er$.

\medskip\textbf{Relations between angles and distances}\qua It is clear that
geodesic rays based at $\er$ and travelling into the core of a
sector travel transversely to walls. We need a quantitative form
of this fact which is the substance of the following:

\begin{lem}\label{3.1} Suppose $\mathfrak{S} \in \what{X}_{\mathfrak{p}}$
 and $\mathfrak{S}\se \ap$ for some apartment $\ap$. Assume that $x\in \mathfrak{S}(\de
)$, $z\in \g _{\mathfrak{S}}$, and $H_z \se \ap$ is a wall
containing $z$. Then $$\measuredangle _z ^\ap (x, H_z ) \geq \sin
^{-1} (\de /2)$$ whenever $d(x,\er )\geq r$ and $d(z,\er) \leq
(\de r)/2$. \end{lem}

 \begin{proof} Notice that $\sa z ^\ap (x,H_z)$ is minimized when $x\in \partial \mathfrak{S}( \de )$,
   $d(x,\er)=r$, and  $H_z$ is parallel to a wall $H_{\er}$ that bounds
 $\mathfrak{S}$.
Therefore, we will assume these three statements are true.
Clearly, $\sa
 z ^\ap (x,H_z)=\sa z ^\ap (x, \pi _{H_z}(x))$ where $\pi _{H_z}\co \ap \rt
H_z$ is the orthogonal projection.
$$d(H_z, H_{\er}) \leq d(z,\er) \leq \frac{\de r}{2},\leqno{\hbox{Note that}} $$
$$d(x,H_{\er}) = d(x,\partial \mathfrak{S}) = \de r.\leqno{\hbox{and}}
$$
Therefore, $$d(x,\pi _{H_z}(x)) = d(x,H_{\er})-d(H_{\er},H_z )
\geq \de r-\frac{ \de r}{2} =\frac{ \de r}{2}.$$ We conclude the
proof by observing that $$ \sa z ^\ap (x,\pi _{H_z}(x)) = \sin
^{-1}\Big[\frac{d(x,\pi _{H_z}(x))}{d(x,z)}\Big]\geq \sin ^{-1}
(\de /2 )
$$ since $d(x,z)\leq d(x,\er) \leq r$.
\end{proof}

The next lemma shows that deep points in the nondegenerate region
of $\xr$ at $\er$ that are separated by a large angle measured at
$\er$ must be a large distance apart. A form of notation we will
use in the proof is $[\er ,z]$ to denote the geodesic segment with
endpoints at $\er $ and $z$.

\begin{lem}\label{3.2} Suppose $x,y \in X_{\mathfrak{p}, \er}(\de )$ and $\Theta
_{\mathfrak{p}}(x,y) \geq 2/( \de r)$, while $d(x,\er) \geq r$ and
$d(y,\er) \geq r$. Then $d(x,y)\geq (\de r)/2$ as long as $\de
\leq 1$. \end{lem}

\proof Choose sectors $\mathfrak{S}_x, \mathfrak{S}_y \in \what
 {X}_{\mathfrak{p}}$ such that $\mathfrak{S}_x \in \Phi _{\mathfrak{p}}(x)$
and $\mathfrak{S}_y \in \Phi _{\mathfrak{p}}(y)$. Let $z \in \xr$
be such that $\g _x \cap \g _y = [\er,z]$. Then, we have $d(\er
,z)\leq ( \de r) /2$ since
$\what{d}_{\mathfrak{p}}(\mathfrak{S}_x, \mathfrak{S}_y)\geq
2/( \de r)$.

 Choose an apartment $\ap _x$ containing $\mathfrak{S}_x$. Note that
$\mathfrak{S}_y \cap \ap _x$ is a convex polyhedron $P$ in $\ap
_x$ that is bounded by walls. Since $z \in
\partial P$, there must be a wall $H_z \se \ap
_x$ such that $z\in H_z$ and $\ap _x - H_z$ has a component which
does not intersect $\mathfrak{S}_y$. Choose a chamber
$\mathfrak{c}_z \se \mathfrak{S}_x$ containing $z$ whose interior
lies in this component, and such that $F=\mathfrak{c}_z \cap \ap
_y$ is a codimension 1 simplex in $\mathfrak{c}_z$.

Let $\mathfrak{c}_y \se \mathfrak{S}_y$ be a chamber containing
$y$. Note that $[z,y] \cup \mathfrak{c}_z \se
{\mathcal{B}}(\mathfrak{c}_z,\mathfrak{c}_y)$, where $
{\mathcal{B}}(\mathfrak{c}_z,\mathfrak{c}_y)$ is the union of
minimal galleries from $\mathfrak{c}_z$ to $\mathfrak{c}_y$.
Hence, $[z,y] \cup \mathfrak{c}_z$ is contained in an apartment
(see e.g. \cite{Br} VI.6). Therefore, $\varrho (\ap _x
,\mathfrak{c}_z) |_{
{\mathcal{B}}(\mathfrak{c}_z,\mathfrak{c}_y)}$ is an isometry,
where $\varrho (\ap _x ,\mathfrak{c}_z) \co \xr \rt \ap _x$ is the
building retraction corresponding to the pair $(\ap _x
,\mathfrak{c}_z)$.

Since $F\se \ap _y$, there is a unique wall $H_z' \se \ap _y$
containing $F$. Since $F \se H_z$ as well, we have $\sa z
^{\ap_y}(y,H_z')=\sa z ^{\ap_x}(\varrho (\ap _x
,\mathfrak{c}_z)(y),H_z)$.

Since $\varrho(\ap _x , \mathfrak{c})$ is distance decreasing, and
since $H_z$ separates $x$ from $\varrho (\ap
_x,\mathfrak{c}_z)(y)$, we have using \fullref{3.1}:
\begin{align*}
d(x,y) & \geq d\big(\varrho (\ap _x ,\mathfrak{c}_z)(x)\,,\,\varrho (\ap _x ,\mathfrak{c}_z)(y)\big) \\
& =d\big(x,\varrho (\ap _x \,,\,\mathfrak{c}_z)(y)\big) \\
& \geq d\big(x,H_z\big) + d\big(\varrho (\ap _x ,\mathfrak{c}_z)(y)\,,\,H_z) \\
& =\sin [\sa z ^{\ap _x} (x,H_z)]d(z,x)+\sin [\sa z ^{\ap _x} (\varrho (\ap _x ,\mathfrak{c}_z)(y),H_z)]d\big(z,\varrho (\ap _x ,\mathfrak{c}_z)(y)\big) \\
& =\sin [\sa z ^{\ap _x} (x,H_z)]d(z,x)+\sin [\sa z ^{\ap _y} (y,H_z')]d(z,y) \\
& \geq  \frac{\de}{2}\Big(d(x,\er)-d(\er,z)\Big)+ \frac{\de}{2}\Big(d(y,\er)-d(\er,z)\Big)\\
& \geq \de r \Big(1-\frac{\de}{2}\Big) \\
& \geq \frac{\de r}{2}.\hskip 3.8in\qed
\end{align*}

Our next lemma states that, after deleting a large compact set, if
the core of two sectors based at $\er$ have a nontrivial
intersection, then the two sectors are close in the boundary metric.

\begin{lem}\label{3.3} Let $\mathfrak{S}_1,\mathfrak{S}_2
 \in \what{X}_{\mathfrak{p}}$, and
suppose that $\mathfrak{S}_1(\de ) \cap \mathfrak{S}_2(\de )\cap
B_{\er}(r)^c \ne \emptyset$. Then
$\what{d}_{\mathfrak{p}}(\mathfrak{S}_1,\mathfrak{S}_2)\leq 2
/(\de r)$. \end{lem}

 \begin{proof} We prove the contrapositive. That is, we assume that
$\g _{\mathfrak{S}_1} \cap \g _{\mathfrak{S}_2}=[\er,z]$ where
$d(\er,z)<(\de r )/2$.

 Choose an apartment $\ap$ with $\mathfrak{S} _2 \se \ap$. We pick
a wall, $H_z$, with $z \in H_z \se \ap$ and such that
$\mathfrak{S}_1 \cap \mathfrak{S}_2 \se \wbar{J}$, where $J$
is a component of $\ap - H_z$ and $\wbar{J}$ is the closure of
$J$.

By \fullref{3.1}, $x \in \mathfrak{S}_2 (\de )\cap B_e(r)^c$ implies
that  $\sa z ^\ap (x,H_z)\geq \sin^{-1}(\de /2)$. Hence, any such
$x$ must be bounded away from $H_z$ and, thus, from
$\wbar{J}$. We have shown  $$\mathfrak{S}_1(\de ) \cap
\mathfrak{S}_2(\de )\cap B_e(r)^c \se \wbar{J} \cap
\mathfrak{S}_2(\de )\cap B_e(r)^c = \emptyset$$ as desired.
\end{proof}

To travel in the nondegenerate space between two deep points
separated by a large angle, one must pass near the origin. More
precisely we have the following:

\begin{lem}[No shifting]\label{3.4} Suppose there is a path $c\co [0,1]
\rt X_e (\de) \cap B_e(r)^c$. Then $\Theta _{\mathfrak{p}} ( c(0),
c(1))\leq 2/(\de r)$. \end{lem}

 \begin{proof} Since $[0,1]$ is compact, it is contained in
finitely many sectors
 $\mathfrak{S}_0,\mathfrak{S}_1,\ldots,\mathfrak{S}_k$ $\in
\what{X}_{\mathfrak{p}}.$ We may assume that these sectors are
ordered so that there exists a partition of $[0,1]$ of the form
$0=t_0<t_1<\ldots<t_{k}=1$ with $c(0) \in \mathfrak{S}_0$, $c(1)\in
\mathfrak{S}_k$, and $c[t_i, t_{i+1}]\se \mathfrak{S}_i$.

 Notice that our partition requires that $c(t_i)\in \mathfrak{S}_i
\cap \mathfrak{S}_{i+1}$. Hence, we can apply \fullref{3.3} to obtain
that $\what
{d}_{\mathfrak{p}}(\mathfrak{S}_i,\mathfrak{S}_{i+1})\leq 2 /(\de
r)$ for all $i$. Therefore, $$ \Theta _{\mathfrak{p}} ( c(0),
c(1))\leq
\what{d}_{\mathfrak{p}}(\mathfrak{S}_0,\mathfrak{S}_k) \leq
{\rm max}
\{\what{d}_{\mathfrak{p}}(\mathfrak{S}_i,\mathfrak{S}_{i+1}) \}
\leq \frac{2}{\de r}.\proved$$
\end{proof}

\textbf{Limit points}\qua Let $\what{X}_\infty$ be the Furstenberg
boundary of $\xy$. That is, we let $\what{X}_\infty$ be the
space of all Weyl chambers up to Hausdorff equivalence. We endow
$\what{X}_\infty$ with the standard metric,
$\what{d}_\infty$, invariant under the stabilizer of $\ey$. We
let $\Phi _\infty \co X_{\infty ,\ey} (\de) \rt
\what{X}_\infty$ be the function that sends a point to its
image at infinity. As $X$ is the product of $\xy$ and $\xr$, we
define $\what{X}= \what{X}_\infty \times
\what{X}_{\mathfrak{p}}$.

A \emph{$\de$--limit point of $\p$ from $e$} is a boundary point
$(\mathfrak{C},\mathfrak{S}) \in \what{X}$, such that there
exists a path $\psi \co [0,\infty) \rt \pn (X_e(\de ))$ that
escapes every compact set, $\lim_{t \to \infty} \Phi_\infty \circ
\p \circ \psi(t)=\mathfrak{C}$, and $\lim_{t \to \infty} \Phi
_{\mathfrak{p}} \circ \p \circ \psi (t)=\{\mathfrak{S} \}$. If this
is the case we call $\psi$ a \emph{limit path} from $e$, and we
write that $\psi$ \emph{limits to} $(\mathfrak{C},\mathfrak{S})$.
We call the set of all $\de$ limit points of $\p$ from $e$, the
$\de$--\emph{limit set of} $\p$ from $e$. We denote the $\de$--limit
set of $\p$ from $e$ by ${\mathcal{L}}_{\p,e} (\de)$.

\medskip\textbf{Existence of nondegenerate visual directions}\qua For the next
result of this section, we return to the material of \fullref{sec2}
and in particular to \fullref{2.9}.

Later we will want to show there are a finite number of limit
points in the limit set of $\p$ to create the finite number of
flats for the conclusion of \fullref{1.2}. This plan will only
succeed if there is a limit point to start with. The results of
\fullref{sec2} were derived for the purpose of showing that limit
points exist. By the Proposition below, we not only know they
exist, we also have precise information on how to construct them.

\begin{prop}[Deep points extended to limit points]\label{3.5}
 Let $\nu _3$ be as in \fullref{2.9}. There is a constant
 $\eta _2=\eta_2(\kappa, \de)$,
such that if $\e \ll \de \ll 1$ and $z \in \pn (X_e(\de ))$ with
$1 \ll _{(\de , \e ,\rho)} r \leq d(z,0)$, then there exists a
boundary point $(\mathfrak{C},\mathfrak{S}) \in
{\mathcal{L}}_{\p,e} (\de/\nu_3)$, such that
$$ \what{d}_{\mathfrak{p}} \big( \mathfrak{S} ,\Phi_{\mathfrak{p}
} \circ \phr (z)\big)\leq \frac{2}{\de r} $$ and $$
\what{d}_\infty \big( \mathfrak{C} ,\Phi_\infty \circ \phy
(z)\big)\leq e^{-\eta _2 r}.$$ \end{prop}

\begin{proof} Let $U$ be the connected component of $\pn (X_e (\de / \nu
_3))\cap B_0 (r)^c$ that contains $z$. From \fullref{2.9} we know
that $U$ is unbounded, so there exists a path $\psi \co [0, \infty
) \rt U$ with $\psi (0) =z $ and such that $\psi$ escapes every
compact set.

Applying \fullref{3.4}, we have that the diameter of $\Phi
_{\mathfrak{p}} \circ \phr \circ \psi ([s, \infty))$ is at most
$2/(\de R_s)$, where $R_s = d (0,\psi ([s,\infty)))$. Notice that
$R_s \to \infty$ as $s \to \infty$, and $$\Phi _{\mathfrak{p}} \circ
\phr \circ \psi \big([t, \infty)\big) \se \Phi _{\mathfrak{p}} \circ
\phr \circ \psi \big([s, \infty)\big)$$ when $ 0\leq s\leq t$.
Therefore, $\lim_{s \to \infty} \Phi _{\mathfrak{p}} \circ \phr
\circ \psi (s)$ exists. Call this limit $\{\mathfrak{S}\}$.

We conclude by remarking that $ \what{d}_{\mathfrak{p}}(
\mathfrak{S} ,\Phi_{\mathfrak{p} } \circ \phr (z))\leq 2 /(\de r)
$ since $$  \Phi _{\mathfrak{p}} \circ \phr \big(z\big) = \Phi
_{\mathfrak{p}} \circ \phr \circ \psi \big(0\big) \in  \Phi
_{\mathfrak{p}} \circ \phr \circ \psi \big([0,\infty )\big)$$ and
$R_0 =r$.

The second part of the proposition is the content of Proposition
5.9 from \cite{E-F1}.
\end{proof}

\textbf{A bound on visual directions for annuli}\qua Once we show
that there is a bound on the number of directions at infinity that
a graded quasiflat can extend in, we can produce a finite
collection of flats that will be our candidates for satisfying the
conclusion of \fullref{1.2}.

Before showing that the number of asymptotic directions a graded
quasiflat travels in is bounded, we will show that the number of
directions is bounded for a quasi-annuli. This bound is
independent of the size of the quasi-annuli. We will then be in a
position to apply the no shifting Lemma in a limiting argument to
show that the same bound exists for the number of directions of a
graded quasiflat.

 Let
$A_R \se \xr$ be the annulus centered at $\er$, with inner radius
$R$ and outer radius $2R$. Let $\phy = \py \circ \p$, and let
$\phr = \pr \circ \p$.

Before proceeding, note that $\pi _\infty
(X_e(\de))=X_{\infty,\ey}(\de)$ and $\pi _{\mathfrak{p}}
(X_e(\de))=X_{\mathfrak{p},\er}(\de)$.

\begin{lem}\label{3.6} The image of $\pr \big[\p(A_R) \cap X_e (\de )\big] $ under
$\Phi _{\mathfrak{p}}$ can be covered by $c_{\mathfrak{p}}=O(1/ \de
^{2m})$ disjoint balls of radius $(4\ka )/(\de ^2 R)$ for $R
> \rho$ and $\e \ll \de $. \end{lem}

 \begin{proof} Let $\mathfrak{S}_i \in \what{X}_{\mathfrak{p}}$ be such that
 $\cup _i B_{\mathfrak{S}_i}(\frac{4\ka}{\de ^2 R}) =
\what{X}_{\mathfrak{p}}$, and $B_{\mathfrak{S}_i}(\frac{4\ka}{\de
^2 R})\cap B_{\mathfrak{S}_j}(\frac{4\ka}{\de ^2 R}) =\emptyset$
if $i \neq j$. That the balls can be chosen to be disjoint is a
consequence of the ultrametric property for
$\what{X}_{\mathfrak{p}}$.

We will twice make use of the fact that if $x \in A_R \cap \p
^{-1} (X_e (\de ) )$, then
\begin{align} d(\phr (x), \er)& = d\big(\p (x) \,,\, (\phy (x),\er)\big) \\
& \geq \de d(\p (x), e) \notag \\
& \geq \frac{\de}{2\ka }D_0(\rho ; x) \notag \\
& \geq \frac{\de R}{2\ka}. \notag
\end{align}
We claim that for any $x \in A_R \cap \p ^{-1} (X_e (\de ) )$, $$
\Phi _{\mathfrak{p}}(\phr (x)) \se B_{\mathfrak{S}_i}\Big(\frac{4\ka}
{\de ^2 R}\Big) \text{ \quad for some }i. $$ Indeed, suppose
$\mathfrak{Z},\mathfrak{Y} \in \Phi _{\mathfrak{p}}(\phr (x))$, and that
$\mathfrak{Z} \in B_{\mathfrak{S}_i}(\frac{4\ka}{\de ^2 R})$. Notice
that $\phr (x) \in X_{\mathfrak{p}, \er} (\de)$, so we can apply (1)
and \fullref{3.3} to obtain
$$\what{d} _{\mathfrak{p}} (\mathfrak{Z},\mathfrak{Y}) \leq \frac{4\ka }{\de ^2 R}.$$
$$\what{d}_{\mathfrak{p}} (\mathfrak{Y},\mathfrak{S}_i) \leq \max \{
\what{d}_{\mathfrak{p}} (\mathfrak{Y},\mathfrak{Z}),
\what{d}_{\mathfrak{p}} (\mathfrak{Z},\mathfrak{S}_i )\} \leq
\frac{4\kappa }{\de ^2 R}\leqno{\hbox{Therefore,}}
$$ as claimed.

Suppose $i\neq j$. If $\Phi _{\mathfrak{p}}(\phr (x))\se
B_{\mathfrak{S}_i}(\frac{4\ka}{\de ^2 R})$ and $\Phi
_{\mathfrak{p}}(\phr (y))\se B_{\mathfrak{S}_j}(\frac{4\ka}{\de ^2
R})$ for a pair of points $x,y \in A_R \cap \p ^{-1} (X_e (\de )
)$, then $B_{\mathfrak{S}_i}(\frac{4\ka}{\de ^2 R}) \cap
B_{\mathfrak{S}_j}(\frac{4\ka}{\de ^2 R})= \emptyset$. Hence, by
the ultrametric property of $\what{X}_{\mathfrak{p}}$ we have
$$ \what{d} _{\mathfrak{p}} \big( \Phi _{\mathfrak{p}} \circ \phr (x)
\,,\, \Phi _{\mathfrak{p}} \circ \phr (y)\big) \geq \frac{4 \ka
}{\de ^2 R}=\frac{2}{\de (\de R/2\ka )}.$$ $$d\big(\phr
(x)\,,\,\phr (y)\big) \geq \frac{\de (\de R/2\ka )}{2}=\frac{\de
^2 R}{4 \ka}\leqno{\hbox{Therefore,}} $$ by (1) and \fullref{3.3}. Thus,
\begin{align*} d(x,y) & \geq \frac{1}{ \ka} d(\p(x), \p
(y))-\e D_0(\rho;x,y) \\
 & \geq \frac{1}{ \ka} d(\phr (x), \phr
(y))-\e D_0(\rho;x,y) \\
 & \geq \frac{\de ^2 R}{4 \ka ^2} - \e 2R \\
 & \geq \frac{\de ^2 R}{5 \ka ^2}.
\end{align*} In summary, we have shown that \begin{align}
d( B_i , B_j) \geq \frac {\de ^2 R}{5 \kappa ^2} \quad \quad (i
\neq j)
\end{align} where $$B_i = A_R \cap \pn \Big[ \pr ^{-1} \Big[\Phi _{\mathfrak{p}}^{-1}
\Big [B_{\mathfrak{S}_i}\Big(\frac{4\ka}{\de ^2 R}\Big)\Big]\Big]
\cap X_e (\de ) \Big].$$
If $\mu _m $ is Lebesgue measure on $\mathbb{E}^m$, then
\begin{align} \mu _m\big[A_R \cap \pn (X_e (\de
))\big] \leq \mu_m \big[A_R \big] < \mu_m \big[B_0(1)\big] (2R)^m.
\end{align} Combining (2) and (3) tells us that the number of
nonempty $B_i$ is bounded above by $$\frac{(10 \kappa ^2)^m
(2R)^m}{(\de ^2 R)^m}=\frac{20^m \kappa ^{2m}}{ \de ^{2m}}.\proved$$
\end{proof}

We will also need to know that projecting onto the symmetric space
factor will produce a bound on the visual angles there. This is
Lemma 4.2 in \cite{E-F1} which we state as

\begin{lem}\label{3.7} 
There exists a constant $\eta_3=\eta_3(\kappa, \de )$, such
that the image of\break $\py \big[\p(A_R) \cap X_e (\de )\big] $ under
$\Phi _\infty$ can be covered by $c_\infty=O(1/ \de ^{2m})$ balls
of radius $e^{- \eta_3 R}$ for $1 \ll _{(\rho ,\de)} R$ and $\e
\ll \de $. \end{lem}

Note that in \cite{E-F1} there is no building factor. Thus, the
statement of Lemma 4.2 in \cite{E-F1} does not mention the
projection map $\py$. Also note that the number of balls in
\cite{E-F1} Lemma 4.2 is bounded by the smaller term $O(1/\de
^m)$. When projecting, a factor of $\de$ makes its way into the
proof from the inequality $d(\py (x), \ey)\geq \de d(x, e)$ for $x
\in X_e(\de )$. The extra factor of $\de$ influences $c_\infty$ by
adjusting the bound from $O(1/\de ^m)$ to $O(1/ \de ^{2m})$, and
our constant $\eta _3$ is proportional to the corresponding
constant in \cite{E-F1}. Aside from these minor adjustments, the
proof carries through without modification.

\medskip\textbf{A bound on visual directions for entire quasiflats}\qua Using
the bound on the number of visual directions for annuli, we are
prepared to pass to the limit and produce a bound for the number
of $\de$--limit points of $\p$.

\begin{prop}[Finite limit set]\label{3.8}  For $\de $ sufficiently small,
$|{\mathcal{L}}_{\p , e}(\de) |<c_\infty c_{\mathfrak{p}}$.
\end{prop}

\begin{proof} Assume there are $c_\infty c_{\mathfrak{p}} +1$ limit points
$\{ ( \mathfrak{C}_i , \mathfrak{S}_i ) \}_{i=1}^{c_\infty
c_{\mathfrak{p}} +1}$. We will arrive at a contradiction.

There are two cases to consider as either $$\big|\{ \mathfrak{C}_i
\}_{i=1}^{c_\infty c_{\mathfrak{p}} +1}\big|>c_\infty \text{\quad or
\quad} \big|\{ \mathfrak{S}_i \}_{i=1}^{c_\infty c_{\mathfrak{p}}
+1}\big|>c_{\mathfrak{p}}.$$ We will begin by assuming the latter.

After possibly re-indexing, let
$\mathfrak{S}_1,\mathfrak{S}_2,\ldots\mathfrak{S}_{c_{\mathfrak{p}+1}}$
be distinct elements of $\{  \mathfrak{S}_i
\}_{i=1}^{c_\infty c_{\mathfrak{p}} +1}$.
Let $\alpha =\min_{i\neq j}\{\what{d}
(\mathfrak{S}_i,\mathfrak{S}_j) \}$. By assumption, there are
paths $$\psi_i \co  [0,\infty) \rt \pn (X_e(\de))$$ such that
$\lim_{t \to \infty} \Phi _{\mathfrak{p}} \circ \phr \circ \psi _i
(t)=\{ \mathfrak{S}_i \}$. Pick $t_i >0$ such that
\begin{align} \bigcup \Phi_{\mathfrak{p}} \circ \phr \circ \psi _i ([t_i,\infty))\se
B_{\mathfrak{S}_i}\Big(\frac{\alpha}{ 2}\Big) \text { \quad for
all } 0 \leq i \leq c_{\mathfrak{p}} +1.
\end{align} We will need a more uniform choice for the $t_i$ to allow us to apply \fullref{3.6}, so we let
$$R = \max \Big\{\frac{8 \kappa } {\alpha \de ^2}, d\big(\psi _1
(t_1),0\big), d\big(\psi _2 (t_2),0\big), \ldots,d\big(\psi
_{c_{\mathfrak{p}+1}} (t_{c_{\mathfrak{p}+1}}),0\big)\Big\}.$$
Then we take $t_i' >0$ such that $d(\psi_i(t_i'),e)=R$ for all $0
\leq i \leq c_{\mathfrak{p}} +1$.

By our choice of $\alpha$,
$$B_{\mathfrak{S}_i}\Big(\frac{\alpha}{2}\Big) \cap
B_{\mathfrak{S}_j}\Big(\frac{\alpha}{2}\Big)=\emptyset \text{
\quad for }i\neq j.$$ Therefore, by (4), $$
B_{\mathfrak{Z}_i}\Big(\frac{\alpha}{2}\Big) \cap
B_{\mathfrak{Z}_j}\Big(\frac{\alpha}{2}\Big)=\emptyset \text{
\quad for }i\neq j, $$ where $ \mathfrak{Z}_i \in
\what{X}_{\mathfrak{p}}$ is a sector containing $\phr \circ \psi
_i (t_i') $. In particular, $\mathfrak{Z}_i \not\in
B_{\mathfrak{Z}_j}(\alpha / 2)$ for $i \neq j$. However, we can
apply \fullref{3.6} to obtain a proper subset $P$ of
$\{1,\ldots,c_{\mathfrak{p}}+1\}$ such that
$$\{\mathfrak{Z}_i\}_{i=1}^{c_{\mathfrak{p}}+1} \se \bigcup_{i\in P}
B_{\mathfrak{Z}_i}\Big(\frac{\alpha}{2}\Big).$$ This is a
contradiction.

If we assume $\big|\{  \mathfrak{C}_i \}_{i=1}^{c_\infty
c_{\mathfrak{p}} +1}\big|>c_\infty$, we can arrive at a similar
contradiction using \fullref{3.7}. The details are carried out in
Proposition 5.2 in \cite{E-F1}.
\end{proof}

\section{Independence of basepoint}\label{sec4}

So far we have limited ourselves by considering a fixed basepoint
$e$. The proof of \fullref{1.2} will require us to hop around from
point to point in our quasiflat with holes, and to treat several
points as basepoints for the nondegenerate space and, hence, for
the limit set of $\p$. We will need to know therefore, that all of
the corresponding nondegenerate spaces and limit sets are
compatible with each other---that they are the same up to minor
modifications of $\de$.

The following lemma is essentially Lemma 5.3 from \cite{E-F1}.

\begin{lem}\label{4.1} Let $r>0$ be given and let $e' \in X$ be such that
$d(e,e')\leq r$. If $x \in \pn (X_e(\de))$ and $d(x,0) \geq
\max\{\rho , (6\kappa r)/\de \}$ for some $x\in \mathbb{E}^m$,
then $x \in$\break $ \pn (X_{e'} (\de /2))$ as long as $\de \leq 1/3$.
\end{lem}

The next lemma is a short technical remark used in the final lemma
of this section.

\begin{lem}\label{4.2} There exists a constant $\nu_4 =\nu_4 (\xr)$ such that if
$\mathfrak{S}\se X_{\mathfrak{p}}$ is a sector based at $e$, and
$\mathfrak{S}' \se \xr$ is a sector based at $e'\in \xr$ with
$\hd(\mathfrak{S},\mathfrak{S}')<\infty$, then there is a sector
$\mathfrak{Z}\se \mathfrak{S}\cap \mathfrak{S}'$ such that
$\hd(\mathfrak{Z},\mathfrak{S})\leq \nu _4 d(e,e')$. \end{lem}

\begin{proof} Let $\mathfrak{S}$ be contained in an apartment $\ap$. Then
there are isometries $a,n_1,n_2,$ $\ldots,n_k \in \text{Isom}(\xr)$ such
that $a$ stabilizes $\ap$, each $n_i$ stabilizes a half-space of
$\ap$ containing a subsector of $\mathfrak{S}'$, and $k$ is
bounded by a constant depending only on $X$.

It is clear that the result holds if $\mathfrak{S}'=a\mathfrak{S}$
or $\mathfrak{S}'=n_i\mathfrak{S}$. Hence the result for the
general $\mathfrak{S}'$ holds by the triangle inequality.
\end{proof}

We are prepared to show that the $\de$--limit set of $\p$ is as
independent of the choice of basepoint as one would expect. First
though we need to identify the boundaries of $\xr$ created using
two different basepoints. Previously we had defined
$\what{X}_{\mathfrak{p}}$ in a way that depended on $\er$. This
was done mostly for notational convenience, but the dependence on
a basepoint would now be a hindrance for us.

Our solution is to give an equivalent definition of
$\what{X}_{\mathfrak{p}}$ as the space of all sectors with
arbitrary basepoints modulo the equivalence that two sectors be
identified if they are a finite Hausdorff distance from each other
(this is equivalent to the condition that the intersection of the
two  sectors contains a third sector). Now the metric on
$\what{X}_{\mathfrak{p}}$ is determined by a choice of a
basepoint (only up to a Lipschitz equivalence though), but the
space $\what{X}_{\mathfrak{p}}$ itself is independent of that
choice.

\begin{lem}\label{4.3} Let $e'=\p (0')$ for some $0'\in \mathbb{E}^m$, and suppose
$\p$ is a $(\kappa , \rho, \e)$ $ \gq$ based at $0'$ as
well as at $0$. If $\de \ll 1$, then $\mathcal{L}_{\p ,e'} (\de)
\se \mathcal{L}_{\p , e} (\de /2)$. \end{lem}

\begin{proof} Suppose $(\mathfrak{C}',\mathfrak{S}') \in \mathcal{L}_{\p
,e'} (\de)$. Then there is a path $\psi \co [0,\infty) \rt$
$ \pn (X_{e'}(\de))$ such that the path $\phr \circ \psi \co
[0,\infty)\rt X_{\pi _{\mathfrak{p}}(e')}(\de)$ escapes every
compact set and limits to $\{\mathfrak{S}'\}$ when observed from
$\pi _{\mathfrak{p}} (e')$.

Let $\mathfrak{S}$ be the sector based at $\er$ such that
$\hd(\mathfrak{S}',\mathfrak{S})<\infty$. Our goal is to show that
$\phr \circ \psi$ limits to $\mathfrak{S}$ when observed from
$\er$.

To this end, for a given $t>0$, let $\mathfrak{S}_t$ be a sector
based at $\er$ such that $\phr \circ \psi (t) \in \mathfrak{S}_t$.
Let $\mathfrak{S}_t'$ be a sector based at $\pi _{\mathfrak{p}}(e')$
such that $\hd(\mathfrak{S}_t',\mathfrak{S}_t)<\infty$. Note that,
by \fullref{4.2}, $\phr \circ \psi (t) \in \mathfrak{S}_t'$ for
sufficiently large values of $t$. Hence, the family
$\mathfrak{S}_t'$ limits to $\mathfrak{S}'$ from the vantage point
of $\pi _{\mathfrak{p}} ( e')$.

Therefore, for any number $r>0$ and sufficiently large values of
$t$, we have $\g _{\mathfrak{S}'}(r) \in \mathfrak{S}_t'$. Recall
that $\g _{\mathfrak{S}'}$ is the geodesic ray in $\mathfrak{S}'$
based at $\pi _{\mathfrak{p}}(e ')$ that travels down the center of
$\mathfrak{S}'$ and is used for measuring distances between points
in $\what{X}_{\mathfrak{p}}$ from the vantage point of $\pi
_{\mathfrak{p}} (e')$.

By Lemmas \ref{4.1} and \ref{4.2}, $\g _{\mathfrak{S}'}(r) \in
\mathfrak{S}_t(\de /2) \cap \mathfrak{S}(\de /2)$. Now applying
the no shifting Lemma gives us that
$$\what{d}_{\mathfrak{p}}(\mathfrak{S}_t,\mathfrak{S})\to 0$$ as
$t \to \infty$. Therefore, $$\lim_{t \to \infty} \Phi_{\mathfrak{p}}
\circ \phr \circ \psi (t)=\{\mathfrak{S}\}$$ as desired.

For the symmetric space part of the proof, see Lemma 5.4 of
\cite{E-F1}.
\end{proof}

\section[Proof of Theorem 1.2]{Proof of \fullref{1.2}}\label{sec5}

Using the tools we have assembled thus far (in particular
large-scale homology of pinched sets, the no shifting Lemma,
extending deep points to limit points, the bound on limit points,
and the independence of basepoints) we can retrace the proof of
Eskin--Farb given in \cite{E-F1} to prove the quasiflats with
holes theorem. Since this proof is essentially contained in
\cite{E-F1}, we will at times only sketch the arguments.

\begin{proof}[Proof of \fullref{1.2}] Since $\Omega _{(\e  ,
\rho ')} \se \Omega _{(\e , \rho)}$ when $\rho ' < \rho$, we may
assume that $1 \ll _{(C)} \rho $. We let $\e $ and $\de$ be
positive numbers such that
  $ \e  \ll \de \ll 1$.

 As in the proof of Theorem 8.1 of \cite{E-F1}, if $x \in
\Omega_{(\e  ,\rho )}$, we can use a connect-the-dots construction
to define a continuous map $\p _x \co \mathbb{E}^m \rt X$ such
that $d(\p _x (y), \varphi (y)) \leq O(\e  )D_x(\rho ;y)$. Hence,
$\p_x$ is a $(\kappa, \rho, O(\e ) )$ $\gq$ based at $x$.

Let $\partial X$ be the Tits building for $X$. Because
$\what{X}$ can be identified with the simplices of maximal
dimension in $\partial X$, we can measure their distances under
the Tits metric. It is well known that if a pair of points in
$\what{X}$ have maximal Tits distance (``opposite points"),
then there is a unique flat that contains the pair up to Hausdorff
equivalence. Let $F_1,\ldots,F_M$ be the flats so obtained from pairs
of opposite points in $\mathcal{L}_{\p_x, \p_x(x)}(\de)$. Note
that $M \leq (c_\infty c_{\mathfrak{p}})^2$ where $c_\infty$ and
$c_{\mathfrak{p}}$ are as in Lemmas \ref{3.6} and \ref{3.7}.

We will show that $\p_x(x)$ is contained in a bounded neighborhood
of $\cup _{i=1}^M F_i$, but first we want to demonstrate that the
limit set, and hence our choice of flats, is independent of $x$.

Suppose $z\in \Omega _{(\e , \rho )}$ and $\p _z$ is constructed
as $\p_x$ to be a $(\kappa ,  \rho ,O(\e ))$ $\gq$ of
$\mathbb{E}^m$ based at $z$.

By construction, we have for any point $y \in X$:
\begin{align} d(\p_z(y),\p_x(y))\leq O(\e )\big(D_z(\rho ;y)+D_x(\rho
;y)\big). \end{align} It follows that $\p _z$ is a $(2 \kappa ,
\rho +2d(x,z), O(\e ))$ $\gq$ based at $x$. Hence, we obtain
through \fullref{4.3} that
$$\mathcal{L}_{\p_z, \p_z(z)}(\de) \se \mathcal{L}_{\p_z,
\p_z(x)}(\de /2).$$
If $(\mathfrak{C},\mathfrak{S}) \in \mathcal{L}_{\p_z,
\p_z(x)}(\de /2)$, then there is a corresponding limit path $\psi
\co [0,\infty) \rt \pn _z (X_{\p_z(x)}(\de /2))$ that limits to
$(\mathfrak{C},\mathfrak{S})$. 

It follows from (5) that 
$\psi (t) \in \pn _x(X_{\p_x(x)} (\de /4))$ for sufficiently large
values of $t$.

By projecting $\psi$ onto factors and applying \fullref{3.2} of this
paper and Lemma 4.1.i of \cite{E-F1} respectively, we see that
$$\Theta_{\mathfrak{p},\p _x(x)}\Big( \pi
_{\mathfrak{p}} \circ \p _x \circ \psi (t),  \pi _{\mathfrak{p}}
\circ \p _{z} \circ \psi (t)\Big) \to 0$$ and $$\Theta_{\infty ,
\p _x(x)}\Big( \py \circ \p _x \circ \psi (t), \py \circ \p _z
\circ \psi (t)\Big) \to 0$$ as $t \to \infty$. The function
$\Theta_{\mathfrak{p},\p _x(x)}$ above is defined analogously to
$\Theta_{\mathfrak{p}}$ with a basepoint of $\pi _{\mathfrak{p}}
(\p_x(x))$ rather than $\er$, and $\Theta _{\infty , \p_x (x)}$ is
the Furstenberg angle between points in $\xy$ measured at the
point $\pi _\infty(\p_x(x))$.

Therefore, it must be that $(\mathfrak{C},\mathfrak{S}) \in
\mathcal{L}_{\p_x, \p_x(x)}(\de /4)$, and hence,
$$\mathcal{L}_{\p_z, \p_z(x)}(\de /2) \se \mathcal{L}_{\p_x,
\p_x(x)}(\de /4).$$
Joining this inclusion with the previous inclusion of limit sets
we have $$\mathcal{L}_{\p_z, \p_z(z)}(\de) \se \mathcal{L}_{\p_x,
\p_x(x)}(\de /4).$$ This shows that our choice of flats is well
defined up to replacing $\de$ with $\de /4$.

Now we return to the task of showing that $\p _x(x)$ is within a
bounded distance of $\cup _{i=1} ^M F_i$. For the remainder of the
proof we let $\p =\p _x$.

For a fixed $1 \ll _{(\rho,\de,\e)}  R$ there must be a $y \in \pn
(X_{\p(x)}(\de))$ such that $d(\p(y),\p(x))=2R$.
Otherwise, we could apply \fullref{2.8} with $W\se \mathbb{E}^m$ equal
to the sphere centered at $x$ with radius $d(x,y)$ to obtain a
contradiction.

Let $e=(\ey , \er)$ be the midpoint of the geodesic between $\p
(y)$ and $\p(x)$ so that $\p (x), \p(y) \in X_e (\de)$. We project
to each factor. Again we will examine the case of a building.

By \fullref{3.5}, there are limit points
$(\mathfrak{C}_i,\mathfrak{S}_i)\in \mathcal{L}_{\p, e}(\de)$ for
$i=1,2$ such that $\what{d}_{\mathfrak{p}}(\mathfrak{S}_1,\Phi
_{\mathfrak{p}} \circ \phr (x))\leq 2/ (\de R)$ and
$\what{d}_{\mathfrak{p}}(\mathfrak{S}_2,\Phi _{\mathfrak{p}} \circ
\phr (y))\leq 2/ (\de R)$. This implies that in the link at
$\er$---denoted by $L_{\er} \se \xr$---the chambers
$\mathfrak{S}_1 \cap L_{\er}$ and $\mathfrak{S}_2 \cap L_{\er}$
are opposite. Therefore, $\mathfrak{S}_1$ and $\mathfrak{S}_2$ are
opposite in $\what{X}_{\mathfrak{p}}$ under the Tits metric, and
there is a unique apartment $\ap ^{12} \se \xr$ that contains
subsectors of $\mathfrak{S}_1$ and $\mathfrak{S}_2$.

We also note that the geodesic segments $[\er ,\phr (x)]$ and
$[\er,\phr (y)]$ can be extended to geodesic rays $\g _x \se
\mathfrak{S}_1$ and $\g _y \se \mathfrak{S}_2$ respectively. The
bi-infinite path $\g _x \cup \g _y$ is a local geodesic, so it is
a global geodesic which we name $\g$.

As $\g$ is a convex subset of Euclidean space, it is contained in
an apartment $\ap ' \se \xr$. Since $\g \se \ap'$, we have that
$\ap'$ contains subsectors of $\mathfrak{S}_1$ and
$\mathfrak{S}_2$. Hence,
$$\er \in \g \se \ap ' =\ap ^{12}.$$ Therefore, $$d\big( \phr
(x),\ap^{12}\big) \leq d\big(\phr (x) ,\er\big) \leq R.$$
In the proof of Theorem 1.1 of \cite{E-F1}, it is shown that
there is a constant $\Lambda $, depending only on $\xy$, and a
flat $F^{12} \se \xy$ that contains $\mathfrak{C}_1$ and
$\mathfrak{C}_2$ up to Hausdorff equivalence, and such that
$$d\big(\phy (x), F^{12}\big) \leq \frac{1}{2}\big(\kappa R +C\big) + \Lambda.$$
Combining this inequality with its building analogue above yields:
$$d\big(\p (x),F^{12} \times \ap^{12}\big) \leq \sqrt{R^2 +
\Big(\frac{1}{2}\big(\kappa R +C\big) + \Lambda \Big)^2}.$$
The proof of \fullref{1.2} is completed by observing that $F^{12}
\times \ap^{12} \se X$ is the unique flat that contains
$(\mathfrak{C}_1,\mathfrak{S}_1)$ and
$(\mathfrak{C}_2,\mathfrak{S}_2)$ up to Hausdorff equivalence.
Hence, $F^{12} \times \ap^{12} \in \{F_i\}_{i=1}^M$. We take the
constant $N$ in the statement of \fullref{1.2} to be $\sqrt{R^2 +
(1/2(\kappa R +C\big) + \Lambda )^2}$.
\end{proof}

\bibliographystyle{gtart} \bibliography{link}

\end{document}